\documentclass[12pt]{article}

\usepackage{amsmath, epsfig, cite}
\usepackage{amssymb}
\usepackage{amsfonts}
\usepackage{latexsym}
\usepackage{float}
\usepackage{color}
\usepackage{booktabs}
\usepackage{graphicx}
\usepackage{caption}

\newtheorem{thm}{Theorem}[section]

\newtheorem{conj}[thm]{Conjecture}
\newtheorem{lem}[thm]{Lemma}

\numberwithin{equation}{section}

\newcommand{\qed}{{\hfill$\square$}\medskip}
\UseRawInputEncoding
\begin{document}

\begin{center}
{\Large\bf Supercongruences involving Ap\'ery-like numbers\\[7pt]
  and Bernoulli numbers}
\end{center}

\vskip 2mm \centerline{Ji-Cai Liu}
\begin{center}
{\footnotesize Department of Mathematics, Wenzhou University, Wenzhou 325035, PR China\\[5pt]
{\tt jcliu2016@gmail.com} \\[10pt]
}
\end{center}

\vskip 0.7cm \noindent{\bf Abstract.}
We establish supercongruences for two kinds of Ap\'ery-like numbers, which involve Bernoulli numbers and Bernoulli polynomials. Conjectural supercongruences of the same type for another four kinds of Ap\'ery-like numbers are also proposed.


\vskip 2mm
\noindent{\it MR Subject Classifications}: 11A07, 11B65, 11M41

\section{Introduction}
In 1979, Ap\'ery \cite{apery-asterisque-1979} introduced the two sequences $\{a_n\}_{n\ge 0}$ and $\{b_n\}_{n\ge 0}$ through $3$-term recurrences in the proof of the irrationality of $\zeta(3)$ and $\zeta(2)$:
\begin{align}
&(n+1)^3 a_{n+1}- (2n+1)(17n^2+17n+5)a_n+n^3 a_{n-1}=0,\quad (a_0=1,a_1=5),\label{a-1}\\[5pt]
&(n+1)^2b_{n+1}-(11n^2+11n+3)b_n-n^2b_{n-1}=0,\quad(b_0=1,b_1=3).\label{a-2}
\end{align}
The two sequences $\{a_n\}_{n\ge 0}$ and $\{b_n\}_{n\ge 0}$ are known as the famous Ap\'ery numbers, which possess the binomial sum formulae:
\begin{align*}
&a_n=\sum_{k=0}^n{n\choose k}^2{n+k\choose k}^2,\\[5pt]
&b_n=\sum_{k=0}^n{n\choose k}^2{n+k\choose k}.
\end{align*}

Zagier \cite{zagier-b-2009} investigated the following recurrence related to \eqref{a-2}:
\begin{align}
(n+1)^2u_{n+1}-(An^2+An+\lambda)u_n+Bn^2u_{n-1}=0,\quad (u_{-1}=0,u_0=1),\label{rec-apery-1}
\end{align}
and searched for triples $(A,B,\lambda)\in \mathbb{Z}^3$ such that the solution of the recurrence \eqref{rec-apery-1} is an integer sequence $\{u_n\}_{n\ge 0}$. Six sporadic sequences are found in Zagier's search, which include the desired solution $\{b_n\}_{n\ge 0}$ and the following sequence (named ${\bf C}$):
\begin{align*}
C^{*}_n=\sum_{k=0}^n{n\choose k}^2{2k\choose k},
\end{align*}
for non-negative integers $n$.

Almkvist and Zudilin \cite{az-b-2006} studied the other recurrence related to \eqref{a-1}:
\begin{align}
(n+1)^3u_{n+1}-(2n+1)(an^2+an+b)u_n+cn^3u_{n-1}=0,\quad (u_{-1}=0,u_0=1),\label{a-3}
\end{align}
and searched for triples $(a,b,c)\in \mathbb{Z}^3$ such that the solution of the recurrence \eqref{a-3}
is an integer sequence $\{u_n\}_{n\ge 0}$. They also found six sporadic sequences, which include the desired solution $\{a_n\}_{n\ge 0}$ and the Domb numbers (named $(\alpha)$):
\begin{align*}
D_n=\sum_{k=0}^n{n\choose k}^2{2k\choose k}{2n-2k\choose n-k},
\end{align*}
for non-negative integers $n$.

Let $\mathbb{N}$ denote the set of non-negative integers and $\mathbb{Z}^{+}$ denote the set of positive integers. Let $p\ge 5$ be a prime and $n,m\in \mathbb{Z}^{+}$.
Chan et al. \cite{ccs-ijnt-2010} showed that
\begin{align*}
D_{np}\equiv D_n\pmod{p^3}.
\end{align*}
Osburn and Sahu \cite{os-aam-2011,os-facm-2013} proved that
\begin{align*}
D_{np^m}\equiv D_{np^{m-1}}\pmod{p^{3m}},
\end{align*}
and
\begin{align*}
C^{*}_{np^m}\equiv C^{*}_{np^{m-1}}\pmod{p^{2m}}.
\end{align*}

The Bernoulli numbers $B_n$ and the Bernoulli polynomials $B_n(x)$ ($n\in \mathbb{N}$) are defined by the generating functions:
\begin{align*}
&\frac{z}{e^z-1}=\sum_{n=0}^{\infty}B_n \frac{z^n}{n!},\\[5pt]
&\frac{ze^{xz}}{e^z-1}=\sum_{n=0}^{\infty}B_n(x)\frac{z^n}{n!}.
\end{align*}

In 2020, Sun \cite[Conjectures 5.1 and 5.3]{sunzh-a-2020} conjectured some supercongruences involving Bernoulli numbers and Bernoulli polynomials for the Ap\'ery-like numbers $\{D_n\}_{n\ge 0}$ and $\{C^{*}_{n}\}_{n\ge 0}$.
The conjectural supercongruence concerning $C^{*}_p$ was proved by Mao \cite{mao-ijnt-2017}, which is earlier than Sun's conjecture.
The conjectural supercongruences concerning $D_p$ were proved by Zhang \cite{zhang-rmjm-2022}.
The conjectural supercongruences concerning $C^{*}_{np}$ and $D_{np}$ for $n=2,3$ were proved by Mao and Wang \cite{mw-pems-2024}.

For $n\in\mathbb{N}$ and $r,s\in \mathbb{Z}^{+}$, let
\begin{align*}
D_n^{(r,s)}=\sum_{k=0}^n{n\choose k}^r\left( {2k\choose k}{2n-2k\choose n-k}\right)^s.
\end{align*}
Note that $D_n^{(2,1)}$ coincides with $D_n$.

The motivation of the paper is to establish supercongruences satisfied by
$D_{np}^{(r,s)}$ and $C^{*}_{np}$ for primes $p\ge 5$ and $r,s,n\in \mathbb{Z}^{+}$ with $r\ge 2$,
which involve Bernoulli numbers and Bernoulli polynomials.

The rest of the paper is organized as follows. We state the main results in the next section.
To prove the main results, we need some preliminary results which are established in Section 3. The proofs of the main results are given in Sections 4 and 5. In the final section, we propose two conjectural supercongruences concerning another four kinds of Ap\'ery-like numbers (named $\bf B,\bf F,(\delta),(\zeta)$).

\section{Main results}
\begin{thm}\label{t-2}
Let $p\ge 5$ be a prime and $r,s,n\in \mathbb{Z}^{+}$ with $r\ge 2$. Then
\begin{align*}
D_{np}^{(r,s)}\equiv D_n^{(r,s)}+p^3B_{p-3} \mathcal{D}^{(r,s)}_n\pmod{p^4},
\end{align*}
where $\mathcal{D}^{(r,s)}_n$, independent of $p$, are given by
\begin{align*}
\mathcal{D}^{(2,1)}_n&=8\sum_{k=0}^{n-1}{n\choose k}^2 {2k\choose k}{2n-2k-2\choose n-k-1}n(n-k)^2\\[5pt]
&-\frac{1}{3}\sum_{k=0}^n{n\choose k}^2{2k\choose k}{2n-2k\choose n-k}\left(2nk(n-k)+2k^3+2(n-k)^3\right),
\end{align*}
and
\begin{align*}
&\mathcal{D}^{(r,s)}_n\\[5pt]
&=-\frac{1}{3}\sum_{k=0}^n{n\choose k}^r\left({2k\choose k}{2n-2k\choose n-k}\right)^s\left(rnk(n-k)+2sk^3+2s(n-k)^3\right)\quad\text{for $r+s\ge 4$.}
\end{align*}
\end{thm}

Let $(\frac{a}{p})$ denote the Legendre symbol for any integer $a$ and odd prime $p$.
\begin{thm}\label{t-3}
Let $p\ge 5$ be a prime and $n\in \mathbb{Z}^{+}$. Then
\begin{align*}
C_{np}^{*}\equiv C_n^{*}+p^2\left(\frac{p}{3}\right)B_{p-2}\left(\frac{1}{3}\right)\mathcal{C}_n^{*}\pmod{p^3},
\end{align*}
where $\mathcal{C}_n^{*}$, independent of $p$, is given by
\begin{align*}
\mathcal{C}_n^{*}=\frac{1}{2}\sum_{k=0}^{n-1}{n\choose k}^2{2k\choose k}(n-k)^2.
\end{align*}
\end{thm}

\section{Preliminaries}
\begin{lem} (\cite[Remark 4 (1)]{ht-jnt-2008})
Let $p\ge 5$ be a prime and $n,k\in \mathbb{N}$ with $n\ge k$. Then
\begin{align}
{np\choose kp}\equiv {n\choose k}\left(1-\frac{1}{3}nk(n-k)p^3B_{p-3}\right)\pmod{p^4}.\label{wolstenholme}
\end{align}
\end{lem}

\begin{lem} (\cite[(33)]{mt-jnt-2018})
For any prime $p\ge 5$, we have
\begin{align}
\sum_{j=1}^{p-1}\frac{1}{j^2}{2j\choose j}\equiv \frac{1}{2}\left(\frac{p}{3}\right)B_{p-2}\left(\frac{1}{3}\right)\pmod{p}.\label{cc-1}
\end{align}
\end{lem}

\begin{lem}
For any prime $p\ge 5$, we have
\begin{align}
&\sum_{j=1}^{p-1}\frac{1}{j^2} \equiv \frac{2}{3}pB_{p-3}\pmod{p^2},\quad\quad &\sum_{j=1}^{p-1}\frac{1}{j^3}\equiv 0\pmod{p}, \label{new-har-6}\\[5pt]
&\sum_{j=1}^{(p-1)/2}\frac{1}{j^2}\equiv \frac{7}{3}pB_{p-3}\pmod{p^2},\quad\quad
&\sum_{j=1}^{(p-1)/2}\frac{1}{j^3}\equiv -2B_{p-3}\pmod{p},\label{new-har-7}
\end{align}
where $H_n$ denotes the $n$th harmonic number $H_n=\sum_{k=1}^n\frac{1}{k}$.
\end{lem}
Note that \eqref{new-har-6} and \eqref{new-har-7} follow from \cite[page 353]{lehmer-am-1938} and \cite[Corollary 5.2]{sunzh-dam-2000}, respectively.

\begin{lem}
Let $p\ge 5$ be a prime and $n,k\in \mathbb{N}$ with $n>k$. Then
\begin{align}
&{2kp+2j\choose kp+j}{2np-2kp-2j\choose np-kp-j}\notag\\[7pt]
&\equiv
\begin{cases}
\displaystyle \frac{2p}{j}{2k\choose k}{2n-2k-2\choose n-k-1}\left(2k+1-2n\right)\quad &\text{for $1\le j\le \frac{p-1}{2}$}\\[20pt]
\displaystyle \frac{2p}{j}{2k\choose k}{2n-2k-2\choose n-k-1}\left(2k+1\right)\quad &\text{for $\frac{p+1}{2}\le j\le p-1$}
\end{cases}\label{cc-8}
\pmod{p^2}.
\end{align}
\end{lem}
{\noindent \it Proof.}
Assume that $1\le j\le p-1$.
Recall the following known result \cite[page 19]{liu-bams-2017}:
\begin{align}
{2kp+2j\choose kp+j}\equiv {2k\choose k}{2j\choose j}\left(1+2kp(H_{2j}-H_j)\right)
\pmod{p^2}.\label{cc-9}
\end{align}
By \eqref{cc-9}, we have
\begin{align}
&{2np-2kp-2j\choose np-kp-j}\notag\\[5pt]
&={(2n-2k-2)p+2p-2j\choose (n-k-1)p+p-j}\notag\\[5pt]
&\equiv {2n-2k-2\choose n-k-1}{2p-2j\choose p-j}\left(1+(2n-2k-2)p(H_{2p-2j}-H_{p-j})\right)
\pmod{p^2}.\label{cc-10}
\end{align}
Note that
\begin{align}
{2p-2j\choose p-j}={2p\choose p}{p\choose j}^2 {2j\choose j}^{-1}{2p\choose 2j}^{-1}.\label{cc-11}
\end{align}
Combining \eqref{cc-9}--\eqref{cc-11} gives
\begin{align}
&{2kp+2j\choose kp+j}{2np-2kp-2j\choose np-kp-j}\notag\\[5pt]
&\equiv 2{2k\choose k}{2n-2k-2\choose n-k-1}{p\choose j}^2 {2p\choose 2j}^{-1}\notag\\[5pt]
&\times \left(1+2kp(H_{2j}-H_j)+(2n-2k-2)p(H_{2p-2j}-H_{p-j})\right)\pmod{p^2},\label{cc-12}
\end{align}
where we have used \eqref{wolstenholme}.

For $1\le j\le \frac{p-1}{2}$, we have
\begin{align}
{p\choose j}^2 {2p\choose 2j}^{-1}&=\frac{p(2j)!(p+1)\cdots (2p-2j)}{2j^2(p+1)\cdots (2p-1)}{p-1\choose j-1}^2\notag\\[5pt]
&\equiv \frac{p(2j)!(p+1)\cdots (2p-2j)}{2j^2(p+1)\cdots (2p-1)}\notag\\[5pt]
&\equiv  -\frac{p}{j}\pmod{p^2},\label{cc-13}
\end{align}
where we have used the fact that ${p-1\choose k}\equiv (-1)^k\pmod{p}$ for $1\le k\le p-1$.

In a similar way, we show that for $\frac{p+1}{2}\le j\le p-1$,
\begin{align}
{p\choose j}^2 {2p\choose 2j}^{-1}\equiv \frac{p}{j}\pmod{p^2}.\label{cc-14}
\end{align}

Finally, combining \eqref{cc-12}--\eqref{cc-14}, we arrive at \eqref{cc-8}.
\qed

\section{Proof of Theorem \ref{t-2}}
Note that
\begin{align}
D_{np}^{(r,s)}&=\sum_{k=0}^n{np\choose kp}^r\left( {2kp\choose kp}{2np-2kp\choose np-kp}\right)^s\notag\\[5pt]
&+\sum_{k=0}^{n-1}\sum_{j=1}^{p-1}{np\choose kp+j}^r\left( {2kp+2j\choose kp+j}{2np-2kp-2j)\choose np-kp-j}\right)^s.\label{ee-1}
\end{align}
By \eqref{wolstenholme}, we have
\begin{align}
&\sum_{k=0}^n{np\choose kp}^r\left( {2kp\choose kp}{2np-2kp\choose np-kp}\right)^s\notag\\[5pt]
&\equiv \sum_{k=0}^n{n\choose k}^r\left({2k\choose k}{2n-2k\choose n-k}\right)^s\left(1-\frac{rnk(n-k)+2sk^3+2s(n-k)^3}{3}p^3B_{p-3}\right)\notag\\[5pt]
&=D_{n}^{(r,s)}-\frac{1}{3}p^3B_{p-3}\sum_{k=0}^n{n\choose k}^r\left({2k\choose k}{2n-2k\choose n-k}\right)^s\notag\\[5pt]
&\times\left(rnk(n-k)+2sk^3+2s(n-k)^3\right)\pmod{p^4}.\label{ee-2}
\end{align}
Next, we shall distinguish two cases to determine the double sum modulo $p^4$:
\begin{align*}
\sum_{k=0}^{n-1}\sum_{j=1}^{p-1}{np\choose kp+j}^r\left( {2kp+2j\choose kp+j}{2np-2kp-2j)\choose np-kp-j}\right)^s.
\end{align*}

{\noindent \bf Case 1} $r=2$ and $s=1$.

By \eqref{wolstenholme}, for $1\le j \le p-1$ we have
\begin{align}
{np\choose kp+j}^2&={np\choose kp}^2\frac{((n-k)p-j+1)^2\cdots ((n-k)p-1)^2((n-k)p)^2}{(kp+1)^2\cdots (kp+j)^2}\notag\\[5pt]
&\equiv \frac{(n-k)^2p^2}{j^2}{n\choose k}^2\pmod{p^3}.
\label{ee-3}
\end{align}
Combining \eqref{cc-8} and \eqref{ee-3} gives
\begin{align*}
&\sum_{j=1}^{p-1}{np\choose kp+j}^2 {2kp+2j\choose kp+j}{2np-2kp-2j)\choose np-kp-j}\\[5pt]
&\equiv -4p^3 n(n-k)^2 {n\choose k}^2 {2k\choose k}{2n-2k-2\choose n-k-1}\sum_{j=1}^{(p-1)/2}\frac{1}{j^3}\\[5pt]
&+2p^3\left(2k+1\right)(n-k)^2{n\choose k}^2 {2k\choose k}{2n-2k-2\choose n-k-1}\sum_{j=1}^{p-1}\frac{1}{j^3}\\[5pt]
&\equiv 8p^3 B_{p-3}n(n-k)^2 {n\choose k}^2 {2k\choose k}{2n-2k-2\choose n-k-1}\pmod{p^4},
\end{align*}
where we have used \eqref{new-har-6} and \eqref{new-har-7}. It follows that
\begin{align}
&\sum_{k=0}^{n-1}\sum_{j=1}^{p-1}{np\choose kp+j}^2 {2kp+2j\choose kp+j}{2np-2kp-2j)\choose np-kp-j}\notag\\[5pt]
&\equiv 8p^3 B_{p-3}\sum_{k=0}^{n-1}{n\choose k}^2 {2k\choose k}{2n-2k-2\choose n-k-1}n(n-k)^2 \pmod{p^4}.\label{ee-4}
\end{align}

Finally, combining \eqref{ee-1}, \eqref{ee-2} and \eqref{ee-4}, we complete the proof of the case $r=2$ and $s=1$ for Theorem \ref{t-2}.

{\noindent \bf Case 2} $r+s\ge 4$.

For $1\le j \le p-1$, we have
\begin{align}
{np\choose kp+j}^r&={np\choose kp}^r\frac{((n-k)p-j+1)^r\cdots ((n-k)p-1)^r((n-k)p)^r}{(kp+1)^r\cdots (kp+j)^r}\notag\\[5pt]
&\equiv 0 \pmod{p^r}.\label{ee-5}
\end{align}
By \eqref{cc-8}, we have
\begin{align}
\left( {2kp+2j\choose kp+j}{2np-2kp-2j)\choose np-kp-j}\right)^s\equiv 0\pmod{p^s}\quad
\text{for $1\le j \le p-1$.}\label{ee-6}
\end{align}
It follows from \eqref{ee-5} and \eqref{ee-6} that
\begin{align*}
{np\choose kp+j}^r\left( {2kp+2j\choose kp+j}{2np-2kp-2j)\choose np-kp-j}\right)^s\equiv 0\pmod{p^4},
\end{align*}
and so
\begin{align}
\sum_{k=0}^{n-1}\sum_{j=1}^{p-1}{np\choose kp+j}^r\left( {2kp+2j\choose kp+j}{2np-2kp-2j)\choose np-kp-j}\right)^s\equiv 0\pmod{p^4}.\label{ee-7}
\end{align}

Finally, combining \eqref{ee-1}, \eqref{ee-2} and \eqref{ee-7}, we complete the proof of the case $r+s\ge 4$ for Theorem \ref{t-2}.

\section{Proof of Theorem \ref{t-3}}
Note that
\begin{align}
C_{np}^{*}=\sum_{k=0}^n{np\choose kp}^2{2kp\choose kp}+\sum_{k=0}^{n-1}\sum_{j=1}^{p-1}{np\choose kp+j}^2{2kp+2j\choose kp+j}.\label{ff-1}
\end{align}
By \eqref{wolstenholme}, we have
\begin{align}
\sum_{k=0}^n{np\choose kp}^2{2kp\choose kp}\equiv \sum_{k=0}^n{n\choose k}^2{2k\choose k}=C_n^{*}\pmod{p^3}.
\label{ff-2}
\end{align}

By \eqref{cc-9}, for $1\le j \le p-1$ we have
\begin{align}
{2kp+2j\choose kp+j}\equiv {2k\choose k}{2j\choose j}\pmod{p}.\label{ff-3}
\end{align}
Combining \eqref{ee-3} and \eqref{ff-3} gives
\begin{align*}
&\sum_{j=1}^{p-1}{np\choose kp+j}^2{2kp+2j\choose kp+j}\\[5pt]
&\equiv (n-k)^2p^2{n\choose k}^2{2k\choose k}\sum_{j=1}^{p-1}\frac{1}{j^2}{2j\choose j}\\[5pt]
&\equiv \frac{1}{2}p^2\left(\frac{p}{3}\right)B_{p-2}\left(\frac{1}{3}\right)(n-k)^2{n\choose k}^2{2k\choose k}\pmod{p^3},
\end{align*}
where we have used \eqref{cc-1} in the last step. It follows that
\begin{align}
&\sum_{k=0}^{n-1}\sum_{j=1}^{p-1}{np\choose kp+j}^2{2kp+2j\choose kp+j}\notag\\[5pt]
&\equiv \frac{1}{2}p^2\left(\frac{p}{3}\right)B_{p-2}\left(\frac{1}{3}\right) \sum_{k=0}^{n-1}{n\choose k}^2{2k\choose k}(n-k)^2\pmod{p^3}.\label{ff-5}
\end{align}

Finally, combining \eqref{ff-1}, \eqref{ff-2} and \eqref{ff-5}, we complete the proof of Theorem \ref{t-3}.

\section{Concluding remarks}
The following four sporadic sequences are also found in Zagier's search \cite{zagier-b-2009} and Almkvist--Zudilin's search \cite{az-b-2006}.

\begin{table}[H]
\centering
\scalebox{0.8}{
\begin{tabular}{cc}
\toprule
Name&Formula\\
\midrule
{\bf B}&$u_n=\sum_{k=0}^n(-1)^k3^{n-3k}{n\choose 3k}{3k\choose 2k}{2k\choose k}$ \\[7pt]
{\bf F}&$u_n=\sum_{k=0}^n(-1)^k8^{n-k}{n\choose k}\sum_{j=0}^k{k\choose j}^3$\\[7pt]
$(\delta)$&$u_n=\sum_{k=0}^n(-1)^k3^{n-3k}{n\choose 3k}{n+k\choose k}{3k\choose 2k}{2k\choose k}$ \\[7pt]
$(\zeta)$&$u_n=\sum_{k=0}^n\sum_{l=0}^n{n\choose k}^2{n\choose l}{k\choose l}{k+l\choose n}$ \\
\bottomrule
\end{tabular}
}
\end{table}

It is interesting that, based on numerical calculation, the sequences ${\bf B},{\bf F},(\delta),(\zeta)$ appear to have the supercongruences of the same type, which involve Bernoulli numbers and Bernoulli polynomials.
\begin{conj}
Let $p\ge 5$ be a prime and $\{u_n\}_{n\ge 0}$ be one of the sequences $\bf B$ and $\bf F$.
For all positive integers $n$, we have
\begin{align}
u_{np}\equiv u_n+\frac{1}{2}p^2\left(\frac{p}{3}\right)B_{p-2}\left(\frac{1}{3}\right)\mathcal{U}_n,
\label{conj-1}
\end{align}
where $\{\mathcal{U}_n\}_{n\ge 1}$ is an integer sequence independent of $p$.
\end{conj}

\begin{table}[H]
\caption*{\text{Values of $\{\mathcal{U}_n\}_{n\ge 1}$} for $\bf B$ and $\bf F$}
\centering
\scalebox{0.8}{
\begin{tabular}{cc}
\toprule
Sequence&Values of $\{\mathcal{U}_n\}_{n\ge 1}$\\
\midrule
$\bf B$& $3, 36, 243, 1008, 675,-32076, -355887, -2483136,\cdots$\\[7pt]
$\bf F$& $10, 240, 3780, 49920, 598500, 6752160, 73076640, 767508480,\cdots$ \\
\bottomrule
\end{tabular}
}
\end{table}

\begin{conj}
Let $p\ge 5$ be a prime and $\{u_n\}_{n\ge 0}$ be one of the sequences $(\delta)$ and $(\zeta)$.
For all positive integers $n$, we have
\begin{align}
u_{np}\equiv u_n+\frac{1}{3}p^3B_{p-3}\mathcal{U}_n,
\label{conj-2}
\end{align}
where $\{\mathcal{U}_n\}_{n\ge 1}$ is an integer sequence independent of $p$.
\end{conj}

\begin{table}[H]
\caption*{\text{Values of $\{\mathcal{U}_n\}_{n\ge 1}$} for $(\delta)$ and $(\zeta)$}
\centering
\scalebox{0.8}{
\begin{tabular}{cc}
\toprule
Sequence&Values of $\{\mathcal{U}_n\}_{n\ge 1}$\\
\midrule
$(\delta)$& $18, 432, 4698, 12672,-492750, -10524816, -118670454, -732312576,\cdots$  \\[7pt]
$(\zeta)$& $-4, -288, -11124, -346368, -9625500, -249508512, -6170456124, -147509102592,\cdots$ \\
\bottomrule
\end{tabular}
}
\end{table}

We remark that cases $n=1,2,3$ of \eqref{conj-1} and \eqref{conj-2} were originally conjectured by
Sun \cite[Conjectures 5.1 and 5.3]{sunzh-a-2020}.

\vskip 5mm \noindent{\bf Acknowledgments.}
This work was supported by the National Natural Science Foundation of China (grant 12171370).


\begin{thebibliography}{99}

\small \setlength{\itemsep}{-.8mm}

\bibitem{az-b-2006}G. Almkvist and W. Zudilin, Differential equations, mirror maps and zeta values,
in Mirror symmetry. V, 481--515, AMS/IP Stud. Adv. Math. 38, Amer. Math. Soc., Providence, R.I., 2006.

\bibitem{apery-asterisque-1979}R. Ap\'ery, Irrationalit\'e de $\zeta(2)$ et $\zeta(3)$, Ast\'erisque 61 (1979), 11--13.

\bibitem{ccs-ijnt-2010}H.H. Chan, S. Cooper and F. Sica, Congruences satisfied by Ap\'ery-like numbers,
Int. J. Number Theory 6 (2010), 89--97.

\bibitem{ht-jnt-2008}C. Helou and G. Terjanian, On Wolstenholme's theorem and its converse, J. Number
Theory 128 (2008), 475--499.

\bibitem{lehmer-am-1938}E. Lehmer, On congruences involving Bernoulli numbers and the quotients of Fermat and Wilson, Ann. Math. 39 (1938), 350--360.

\bibitem{liu-bams-2017}J.-C. Liu, Congruences for truncated hypergeometric series ${}_2F_1$, Bull. Aust. Math. Soc. 96 (2017), 14--23.

\bibitem{mao-ijnt-2017}G.-S. Mao, Proof of some congruences conjectured by Z.-W. Sun, Int. J. Number Theory 13 (2017), 1983--1993.

\bibitem{mw-pems-2024}G.-S. Mao and L. Wang, Proof of some conjectural congruences involving Ap\'ery and Ap\'ery-like numbers, Proc. Edinb. Math. Soc., in press, doi: 10.1017/S0013091524000075.

\bibitem{mt-jnt-2018}S. Mattarei and R. Tauraso, From generating series to polynomial congruences, J. Number Theory 182 (2018), 179--205.

\bibitem{os-aam-2011}R. Osburn and B. Sahu, Supercongruences for Ap\'ery-like numbers, Adv. in Appl. Math. 47 (2011), 631--638.

\bibitem{os-facm-2013}R. Osburn and B. Sahu, A supercongruence for generalized Domb numbers, Funct. Approx. Comment. Math. 48 (2013), 29--36.

\bibitem{sunzh-dam-2000}Z.-H. Sun, Congruences concerning Bernoulli numbers and Bernoulli polynomials,
Discrete Appl. Math. 105 (2000), 193--223.

\bibitem{sunzh-a-2020}Z.-H. Sun, Congruences for two types of Ap\'ery-like sequences, preprint (2020), 	arXiv:2005.02081.

\bibitem{zagier-b-2009}D. Zagier, Integral solutions of Ap\'ery-like recurrence equations, in Groups and Symmetries, 349--366, CRM Proc. Lecture Notes 47, Amer. Math. Soc., Providence, R.I., 2009.

\bibitem{zhang-rmjm-2022}Y. Zhang, Some conjectural supercongruences related to Bernoulli and Euler numbers, Rocky Mountain J. Math. 52 (2022), 1105--1126.

\end{thebibliography}
\end{document}